\documentclass[10pt,twosided]{article}
\usepackage[tbtags]{amsmath}
\usepackage{amssymb}
\allowdisplaybreaks[4]
\usepackage{color}

\newcommand{\Abs}[1]{ \Bigl \lvert #1 \Bigr \rvert}

\definecolor{c20}{rgb}{0.,0.7,0.}
\definecolor{c30}{rgb}{0.,0.,1.}
\definecolor{c40}{rgb}{1,0.1,0.7}
\definecolor{c50}{rgb}{1,0,0}
\definecolor{c60}{rgb}{1,0.9,0.1}
\definecolor{c70}{rgb}{0.50,1.00,0.00}

\def\rH#1{{\textcolor{c30}{#1}}}
\def\rH#1{#1}

\def\cZ#1{{\textcolor{c40}{#1}}}
\def\cZ#1{#1}
\def\CC#1{{\textcolor{c40}{#1}}}
\def\CC#1{#1}
\def\ccc#1{{\textcolor{c40}{#1}}}
\def\ccc#1{#1}
\def\cwz#1{{\textcolor{c40}{#1}}}
\def\cwz#1{#1}
\def\wzc#1{{\textcolor{c40}{#1}}}
\def\wzc#1{#1}
\def\wz#1{{\textcolor{c40}{#1}}}
\def\wz#1{#1}

\def\cH#1{\textcolor{c20}{#1}}
\def\cH#1{#1}
\def\eH#1{\textcolor{c20}{#1}}
\def\eH#1{#1}
\def\EH#1{\textcolor{c20}{#1}}
\def\EH#1{#1}

\def\peng#1{{\textcolor{c40}{#1}}}
\def\peng#1{#1}
\def\pzx#1{{\textcolor{c70}{#1}}}
\def\pzx#1{#1}

\def\cE#1{\textcolor{c30}{#1}}
\def\cE#1{#1}

\def\Peng#1{#1}

\usepackage{amssymb,color}

\newcommand{\ve}{\varepsilon}

\newcommand{\ABs}[1]{ \biggl \lvert #1 \biggr \rvert}

\newcommand{\pk}[1]{\mathbb{P} \left\{ #1 \right\} }

\newcommand{\EXP}[1]{\exp \left( #1 \right) }

\newcommand{\R}{\!I\!\!R}

\newcommand{\inr}{\in \R}

\newcommand{\limit}[1]{\lim_{#1 \to   \infty}}

\newcommand{\BQN}{\begin{eqnarray}}
\newcommand{\EQN}{\end{eqnarray}}
\newcommand{\BQNY}{\begin{eqnarray*}}
\newcommand{\EQNY}{\end{eqnarray*}}

\newcommand{\BS}{\begin{sat}}
\newcommand{\ES}{\end{sat}}
\newcommand{\BT}{\begin{theo}}
\newcommand{\ET}{\end{theo}}
\newcommand{\BK}{\begin{korr}}
\newcommand{\EK}{\end{korr}}
\newcommand{\EQD}{\stackrel{d}{=}}

\newcommand{\BD}{\begin{de}}
\newcommand{\ED}{\end{de}}

\newcommand{\BIT}{\begin{itemize}}
\newcommand{\EIT}{\end{itemize}}
\newcommand{\BDI}{\begin{description}}
\newcommand{\EDI}{\end{description}}

\newcommand{\BRM}{\begin{remark}}
\newcommand{\ERM}{\end{remark}}

\newcommand{\BEL}{\begin{lem}}
\newcommand{\EEL}{\end{lem}}

\numberwithin{equation}{section}
\newtheorem{theo}{Theorem}[section]
\newtheorem{sat}[theo]{Proposition}
\newtheorem{de}[theo]{Definition}
\newtheorem{lem}{Lemma}[section]

\newtheorem{korr}[theo]{Corollary}
\newtheorem{remark}[theo]{Remark}

\newcommand{\prooftheo}[1]{ \textsc{Proof of Theorem} \ref{#1} }

\newcommand{\prooflem}[1]{\textsc{Proof of Lemma} \ref{#1}}
\newcommand{\proofkorr}[1]{\textsc{Proof of Corollary} \ref{#1}}

\newcommand{\QED}{\hfill $\Box$}

\def\d{\mathrm{d}}
\newcommand{\EE}[1]{\mathbb{E}\left\{#1\right\}}

\def\I{\operatorname*{\mathbb{I}}}
\newcommand{\nelem}[1]{{Lemma \ref{#1}}}

\newcommand{\netheo}[1]{{Theorem \ref{#1}}}

\newcommand{\nwc}{\newcommand}
\nwc{\COM}[1]{}

\def\IF{\infty}

\begin{document}
\centerline{ \Large \bf Limit \rH{Properties} of Exceedances Point
Processes of }
\centerline{ \Large \bf Scaled Stationary Gaussian Sequences}

\bigskip
\centerline{Enkelejd Hashorva\footnote{University of Lausanne, UNIL-Dorigny, 1015 Lausanne, Switzerland}, Zuoxiang Peng\footnote{School of Mathematics and Statistics, Southwest
University, 400715 Chongqing, China}, Zhichao Weng$^a$   }

\bigskip
\centerline{\today{}}

{\bf Abstract:} \cE{We derive the limiting \rH{distributions} of exceedances
point processes of randomly scaled weakly dependent stationary Gaussian
sequences under some mild asymptotic conditions. \EH{In the literature analogous results are available only} \rH{for}
contracted stationary \rH{Gaussian sequences}. In this paper, we include additionally the case of randomly inflated stationary Gaussian sequences
\rH{with a} Weibullian type random scaling. \ccc{It} turns out that  the maxima and minima of \eH{both contracted and inflated weakly dependent stationary Gaussian \rH{sequences}} are asymptotically independent.}


\section{Introduction}
\eH{Let $X_{n},n\ge 1$ be a standard stationary Gaussian sequence (ssGs) i.e., $X_n$'s are $N(0,1)$ distributed and
$\rho(n)=\EE{X_1X_{n+1}}= \EE{\cwz{X_jX_{n+j}}}$ for any $j\ge 1$.
In \cE{the seminal contribution} \cite{Berman64}, S.M. Berman \rH{proved} that the maxima
$\tilde{M}_{n}=\cE{\max_{1 \le k \le n}} X_{k}$ converges in distribution after normalization to a unit Gumbel random variable,
i.e., }
\BQN\label{Mn}
\lim_{n\to\infty}\pk{\tilde{M}_{n}\le
\tilde{a}_{n}x+\tilde{b}_{n}}=\exp(-\exp(-x))=:\Lambda(x), \quad \forall x\inr,
\EQN
provided \rH{that the} \cE{so-called Berman condition}
\BQN \label{rho}
\limit{n}
\rho(n) \ln n=0 \EQN holds, where the norming constants
$\tilde{a}_{n}$ and $\tilde{b}_{n}$ are given by
\BQNY
\tilde{a}_n=\frac{1}{\sqrt{2 \ln n}}\, \mbox{   and }\,
\tilde{b}_n=\sqrt{2 \ln n} -\frac{\ln   \ln   n + \ln   4 \pi}{2 \sqrt{2 \ln n}}.
\EQNY
Moreover, the
maxima and the minima \eH{$\cwz{\tilde{m}}_n=\min_{1 \le k \le n} X_k$} are asymptotically independent, cf. \cite{Davis79}
and \cite{leadbetter1983extremes}.

In applications, commonly the observations are randomly scaled, say due to some inflation or deflation effects if financial losses are modeled, or
caused by  measurement errors if observations are the outcome of a certain physical \eH{experiment}.
Therefore,
in order to model some general random scaling phenomena applicable to original data,
in this we consider paper $Y=SX, Y_n=S_n X_n,n\ge 1$ assuming that
$S, S_n,n\ge 1$ are independent non-negative random variables with
common distribution function $F$ being further independent of
\cE{the standard Gaussian random variables} $X,X_n,n\ge 1$.

\cE{As shown in \cite{HashorvaWengB} if $F$ has a finite upper endpoint $x_F\in (0,\IF)$ and its survival function is regularly varying, then
the maxima $M_n=\max_{1 \le k\le n} Y_k$ converge in distribution after normalization to a unit Gumbel random variable with distribution function $\Lambda$, provided that the Berman condition holds.\\
If $x_F=\IF$ and $X_n,n\ge 1$ are iid $N(0,1)$ the convergence of maxima $M_n$ is shown under a different normalization in \cite{HashorvaWengA}
assuming further that  $F$ has a \eH{Weibullian} tail behaviour (see below \eqref{eq:SB})}.

\cE{The objective of the paper is twofold: first for $F$ with a Weibullian tail behaviour,
it is of interest to establish the convergence of maxima of \rH{a} randomly scaled  \eH{ssGs} under the Berman condition;  there is no result in the literature covering this case.
Secondly, for both cases $x_F$ is a positive constant, and $x_F=\IF$, we aim at establishing the same result as in \cite{Davis79}, i.e.,
the asymptotic independence of maxima and minima of randomly scaled weakly dependent \rH{ssGs} .\\
Since by using a point process approach also the joint limiting distribution of upper and lower order statistics can be easily established,
we choose in this paper a point process framework considering exceedances point processes. Numerous authors dealt with the asymptotic behavior of
exceedances point processes;  for weakly dependent stationary sequences including Gaussian, see \cite{leadbetter1983extremes, Pit96,  HuPeng99, Faletal2010, arendarczyk2011asymptotics, arendarczyk2011exact, Peng12} and the references therein.}

\cE{For $u_n(s)= a_n s+ b_n,s\inr$, with $a_n>0,b_n\inr$ we shall investigate the weak convergence of bivariate \eH{point} processes
of exceedances of} levels $u_{n}(x)$ and
$-u_{n}(y)$ formed by $Y_n, n\ge1$.
Setting $\xi_{1}(n)=Y_{n},\xi_{2}(n)=-Y_{n}$ for $n\ge1$ we define as in
\cite{Win96} the bivariate exceedances point processes
\begin{equation}\label{eq pp}
\mathbf{N}_n(\mathbf{B},\mathbf{x})=\sum_{d=1}^2\sum_{i=1}^n \I
\left( \xi_{d}(i)
>u_n(x_d), \frac{i}{n} \in B_{d}\right)
\end{equation}
for $\mathbf{B}=\bigcup^2_{d=1}(B_{d}\times \{d\})$ with $B_{d}$ the
Borel set on $(0,1], \Peng{d}=1,2$, where $\I(\cdot)$ denotes the indicator
function. 
The marginal point processes are defined by
\begin{eqnarray*}
N_{n,d}(B_{d},x_d)&=&\sum_{i=1}^n \I \left( \xi_{d}(i) >u_n(x_d), \frac{i}{n} \in B_{d}\right), \quad d=1,2.
\end{eqnarray*}
%
%
%
%
In order to study the weak convergence of $\mathbf{N}_n$ we need to formulate certain assumptions on the random scaling $S$.\\
Our first model concerns the case that $S$ has a Weibullian type tail behaviour with $x_F=\IF$, whereas \rH{the} second one deals with $S$ having a regular tail behaviour at $x_F$. For both cases we investigate the convergence in distribution of $\mathbf{N}_n$, and further, as in \cite{Davis79} we prove that maxima and minima are asymptotically independent.

\COM{ is bounded, it is known that (see \cite{}) under a modified Berman condition the maxima of the randomly scaled stationary
Gaussian sequence still converge to a unit Gumbel random variable with distribution function $\Lambda$. When $S$ has distribution function $F$ with an infinite upper endpoint, large values of $S$ might be too influential meaning that the large values of $Y_n$'s might be caused only by those values. We shall consider therefore only $S$ with very light tail. Our main results below show the weak convergence of $\mathbf{N}_n$ for the two different cases of random scaling. Moreover, we establish as in \cite{Davis79} the asymptotic independence of maxima and minima.}

The rest of the paper is organized as follows. Section 2 gives the main results. Proofs and auxiliary results are displayed in Section 3.

\section{Main Results}
In order to proceed with the main results we need to specify our models for the random scaling $S\ge 0$ with distribution function $F$.
We consider first the case that $S$ has a Weibullian type tail behaviour, i.e.,
 for given positive constants $L, p $
 \BQN
\label{eq:SB} \overline{F}(u)=\pk{S> u}=(1+o(1)) \ccc{g(u)}
\exp(-L u^{p}), \quad u\to \infty,
 \EQN
where $g$ is an ultimately monotone function satisfying $\limit{t}g(tx)/g(t)=x^{\alpha},\forall x>0$ with some $\alpha \inr$.
 Commonly if the latter asymptotic relation holds, then $g$ is referred to as a regularly varying function at infinity with index $\alpha$.
The assumption \eqref{eq:SB} is crucial for finding the tail asymptotics of $Y=S X$, where $S$ and $X$ are independent and $X$ has $N(0,1)$ distribution. Indeed, in view of \cite{arendarczyk2011asymptotics}
\BQN\label{eq:product}
\pk{Y>u}\sim 
(2+p)^{-\frac{1}{2}}g\left(Q^{-1}u^{\frac{2}{2+p}}\right)\EXP{-Tu^{\frac{2p}{2+p}}}, \qquad u\to \IF,
\EQN
where
\BQN\label{eq T Q}
 T:=2^{-1}Q^2+LQ^{-p}, \quad Q:=(Lp)^{1/(2+p)}.
\EQN
Hence \eqref{eq:product} shows that $Y$ has also a Weibullian type distribution. We state next our first result for this Weibullian type scaling model.

\BT\label{th N} Let $X_{n}, n\ge 1$ \cE{be a stationary Gaussian sequence} satisfying \eqref{rho}, and let $\mathbf{N}_n $ be the
bivariate point process given by \eqref{eq pp} with $S_n,n\ge 1$ such that their common distribution function $F$ satisfies
\eqref{eq:SB}. \ccc{\cE{If} further there exist some sequences
$u_n(x),n \ge 1,x\inr$ such that \rH{for any $x\in \mathbb{R}$}
\BQN\label{eq n tail} \lim_{n \to
\IF}n\pk{Y> u_n(x)}= \EXP{-x},
\EQN}
 then $\mathbf{N}_n$ converge in
distribution to a Poisson process $\mathcal{N}$ on
$\bigcap^2_{d=1}((0,1]\times\{d\})$ with intensity
$\mu(\mathbf{B})=\sum_{d=1}^2\exp(-x_d)m(B_d)$, where $m$ denotes the
Lebesgue measure on $\wzc{(}0,1]$. \ET

\BRM If \cE{\eqref{eq:SB} holds with $g(x)= C x^\alpha ,C>0$, then}
in view of \cite{arendarczyk2011asymptotics}
\BQNY \pk{Y>u}\sim
(2+p)^{-\frac{1}{2}}CQ^{-\alpha} u^{\frac{2\alpha}{2+p}}
\EXP{-Tu^{\frac{2p}{2+p}}}, \quad  u\to \IF. \EQNY
Consequently, \eqref{eq n tail} holds according to \pzx{\cite{Embetal1997}} p.155
with $u_n(x)=a_nx+b_n, x\inr$ \eH{and $Q,T$ as in \eqref{eq T Q}}, where 
\BQNY\label{anbn}
a_n=\frac{2+p}{2p}T^{-\frac{2+p}{2p}}(\ln n)^{\frac{2-p}{2p}},
\qquad b_n=\left(\frac{\ln n}{T}\right)^{\frac{2+p}{2p}}
+ \eH{a_n}\left(\frac{\alpha}{p}\eH{\ln ( T^{-1}\ln n)}+\ln (2+p)^{-\frac{1}{2}}CQ^{-\alpha}\right).
\EQNY
\ERM

\eH{Applying} Theorem \ref{th N} we derive \eH{below} the joint \cE{limiting} distribution of
the $k$th maxima and the $l$th minima which are stated as follows.

\BK\label{co M m} For positive integers $k$ and $l$, let $M_n^{(k)}$
and $m_n^{(l)}$ denote the $k$th largest and the $l$th smallest of
$Y_n, n\ge 1$, then under the conditions of Theorem \ref{th N}, for $x,y
\inr$ we have
\BQN \label{corr}
\lim_{n \to \IF} \pk{M_n^{(k)} \le u_n(x), m_n^{(l)} >
-u_n(y)}=\EXP{-\EXP{-x}-\EXP{-y}}\sum_{i=0}^{k-1}\frac{\EXP{-ix}}{i!}\sum_{j=0}^{l-1}\frac{\EXP{-j
y}}{j!}.
\EQN
\EK

Next, we consider the case $S$ has a finite upper endpoint, say $x_F=1$. As in \cite{HashorvaWengB}
 we shall suppose that for any $u \in (\nu, 1)$ with some $\nu \in (0,1)$
 \BQN\label{eq:bound S} \cZ{
\pk{S_{\tau}> u}\geq \pk{S> u} \ge  \pk{S_{\gamma}> u}} \EQN holds
with  $S_{\gamma},S_{\tau}$ two non-negative random variables which
have a regularly varying survival function at 1 with non-negative
index ${\gamma}$ and ${\tau}$, respectively. \eH{By definition} $\ccc{S_{\alpha}}, \eH{\alpha\ge 0}$ is regularly varying at 1 with index $\alpha$ if the distribution function of $S_\alpha$ has upper endpoint equal 1 and further 
{$$ \lim_{u \to \IF}
\pk{ S_{\alpha} > 1- x /u}/ \pk{ S_{\alpha}> 1- 1/u}=x^\alpha, \quad x>0.$$}
The recent contribution \cite{HashorvaWengB} derives the limit distribution of maxima of $Y_i, 1\le i
\le n$ under \eH{the following} modified Berman condition
\BQN\label{eq rho} \lim_{n \to \IF}\rho(n)(\ln
n)^{1+\Delta_{\epsilon}}=0, \EQN
where $\Delta_\epsilon=2(\gamma-\tau)+\epsilon$ and some $\epsilon>0$. \EH{Our last result below} extends the main finding of \cite{HashorvaWengB} establishing  the weak convergence of the bivariate exceedances point process when $S$ is bounded.

\BT\label{th N bound} Let $\mathbf{N}_n$ be defined as in \eqref{eq
pp} with $S_n$ satisfying \eqref{eq:bound S}. If condition \eqref{eq
rho} is satisfied, then $\mathbf{N}_n$ converge in distribution as
$n \to \IF$ to a Poisson process $\mathcal{N}$ on
$\bigcap^2_{d=1}((0,1]\times\{d\})$ with intensity
$\mu(\mathbf{B})=\sum_{d=1}^2\EXP{-x_d}m(B_d)$\Peng{,} where $m$
denotes \cE{the} Lebesgue measure on $\wzc{(}0,1]$. \ET

\BRM
a) Under the assumptions of Theorem \ref{th N bound} for $x,y \inr$ we have that \eqref{corr} holds. Hence in particular the maxima and minima are asymptotically independent in both models for the tail behaviour of $S$.\\
b) If $S$ is regularly varying at 1 with some index $\gamma$, then the claim of \netheo{th N bound} holds under the Berman condition, i.e.,
the modified Berman condition \EH{should be imposed with $\Delta_\ve=0$.}
\ERM

\section{Further Results and Proofs}

\begin{lem}\label{lem20}
\EH{Let $S, Z_n,n\ge 1$ be independent positive random variables \eH{satisfying}
\BQN
\exp(- \widetilde{L}_0 u^{p_1}) \le \pk{S> u} \le \exp(- L_0 u^{p_1}) , \quad
\exp(- \widetilde{L}_n u^{p_2}) \le \pk{Z_n> u} \le \exp(- L_n u^{p_2})
\EQN
for all \wzc{$u$ large} with  $p_1,p_2,\widetilde{L}_n,L_n,n\ge \wzc{0}$ \eH{positive constants} such that $\widetilde{L}_n,L_n \in [a, b], \forall n \ge 0$ with $a <b$ two finite positive constants.}
\wzc{If further $S^*$} is a positive random variable independent of $Z_n,n \ge 1$ satisfying $\lim_{\cwz{u\to \IF}} \pk{S> u}/\pk{S^*> u}=c\in (0,\IF)$,
then \eH{we have} \rH{uniformly in $n$} as $\wzc{u}\to \IF$
\BQNY
 \pk{SZ_n > \wzc{u}} \sim  c\pk{S^* Z_n > \wzc{u}}.
 \EQNY
\end{lem}

\prooflem{lem20}
Let $G_n,n\ge 1$ be the distribution function of $Z_n$.
\eH{By the independence of $S$ and $Z_n$, for} all \wz{$u$}  large 
\begin{eqnarray*}
\overline{H}(u)&:=&\pk{SZ_n>u}\\
&\geq& \pk{S>u^{\frac{p_2}{p_1+p_2}}} \pk{
Z_n>u^{\frac{p_1}{p_1+p_2}}}
\geq
\EXP{-2b u^{\frac{p_1p_2}{p_1+p_2}}}.
\end{eqnarray*}
\Peng{\eH{Further}, for $c_{1}>0$ small enough \eH{and all \wz{$u$} large} we have}
\BQN
\int^{c_1u^{\frac{p_1}{p_1+p_2}}}_0 \pk{S>\frac{u}{s}}\, \d
G_n(s) &\leq & \pk{S>c_1^{-1}u^{\frac{p_2}{p_1+p_2}}} 
 \leq  \EXP{-a c_1^{-p_1}u^{\frac{p_1p_2}{p_1+p_2}}}=o(\overline{H}(\cwz{u}))
\EQN
\Peng{and for some large $c_{2}>0$}
\begin{eqnarray*}
\int_{c_2u^{\frac{p_1}{p_1+p_2}}}^{\infty}
\pk{S>\frac{u}{s}}\,dG_n(s)
\leq\pk{Z_n>c_2u^{\frac{p_1}{p_1+p_2}}}
\leq 
\EXP{-a c_2^{p_2}u^{\frac{p_1p_2}{p_1+p_2}}}
=o(\overline{H}(\cwz{u})).
\end{eqnarray*}
Therefore, for  $\delta_u=c_1
u^{p_1/(p_1+p_2)}, \rH{\lambda_u}=c_2u^{p_1/(p_1+p_2)}$ we have
\BQN \label{sch} \pk{S Z_n> u} \sim \int_{\delta_u}^{\rH{\lambda_u}}
\pk{S> u/s}\, \d G_n(s). \EQN
\eH{Since further $\lim_{\wzc{u}\to \IF} u/\rH{\lambda_u}= \IF$,  for any} $s\in [\delta_u,\rH{\lambda_u}]$ we have $u/s \ge
u/\wz{\lambda_u} \to \IF$ as $u\to \IF$. \eH{Consequently for any $\ve >0$, $s\in [\delta_u,\lambda_u]$}
$$ c(1- \ve)\le \frac{\pk{S> u/s}}{\pk{S^*> u/s}}\le  c(1+ \ve)$$
holds \rH{uniformly in $n$} for all \wzc{$u$} large implying
$$ \pk{S Z_n> u} \sim c \int_{\delta_u}^{\lambda_u} \pk{S^*> u/s}\, \d G_n(s) \sim \pk{S^* Z_n> u} $$
\Peng{as $\wzc{u}\to \IF$ holds also \rH{uniformly in $n$}, and thus the claim} follows. \QED

\BEL\label{add le3.3}
\EH{Let $L_n,n\ge 1$ be as in \nelem{lem20} and let  $Z_n,n\ge 1$ be positive random variables such that}
$$\overline{G}_{n}(z):=\pk{Z_n>z}=\EXP{-\eH{L}_n z^q}$$
\EH{for some  $q>0$ and all $z>0$}. If further \Peng{$Z_n, n\ge 1$ \EH{are} independent \eH{of a  non-negative random variable} $S$
which satisfies} \eqref{eq:SB}, then \wz{we have uniformly in $n$} 
\BQN\label{eq:bound SB}
\pk{SZ_n> \wz{u} } \sim \sqrt{\frac{2\pi Lp}{p+q}}\cwz{A_n^{\frac{p}{2}}}u^{\frac{pq}{2(p+q)}}g\left(\cwz{A_n}u^{\frac{q}{p+q}}\right)
\EXP{-D_nu^{\frac{pq}{p+q}}}
\EQN
\rH{as $\wz{u}\to \IF$,} where $D_n=\left(L+Lpq^{-1}\right)\cwz{A_n^p}$ and $\cwz{A_n}=\left(q\cwz{L}_n\right)^{\frac{1}{p+q}}\left(Lp\right)^{-\frac{1}{p+q}}$.
\EEL

\prooflem{add le3.3}
If \eqref{eq:SB} holds, by Lemma \ref{lem20}, we have for all $\wz{u}$ large
\begin{eqnarray*}
\pk{SZ_n>\rH{u}}&=&\int_{0}^{\infty }\pk{Z_n>\frac{\rH{u}}{s}}\, \d F(s)
\sim
\int_{c_{1}\rH{u}^{\frac{q}{p+q}}}^{c_2\rH{u}^{\frac{q}{p+q}}}\pk{Z_n>\frac{\rH{u}}{s}}\, \d F(s)\\
&\sim&
\CC{\int_{c_{1}\rH{u}^{\frac{q}{p+q}}}^{c_2\rH{u}^{\frac{q}{p+q}}}\EXP{-\cwz{L}_n \rH{u}^q s^{-q}}\, \d F(s)}\\
&\sim&
\CC{\int_{c_{1}\rH{u}^{\frac{q}{p+q}}}^{c_2 \rH{u}^{\frac{q}{p+q}}}\EXP{-\cwz{L}_n \rH{u}^q s^{-q}}\, \d (g(s)\EXP{-Ls^p})}.
\end{eqnarray*}
\cE{Using similar arguments as in the proof of} Theorem 2.1 in \cite{HashorvaWengA} we obtain as \rH{$u\to \IF$}
\begin{eqnarray*}
\pk{SZ_n>\rH{u}} &\sim& Lp\int_{c_{1}\rH{u}^{\frac{q}{p+q}}}^{c_2\rH{u}^{\frac{q}{p+q}}}
s^{p-1}g(s)\EXP{-\cwz{L}_n\rH{u}^qs^{-q}-Ls^p}\, \d s\\
&=&Lp \cwz{A_n^p} \rH{u}^{\frac{qp}{p+q}}
\int_{c_1A_n}^{c_2A_n}
z^{p-1}g\left(A_n \rH{u}^{\frac{q}{p+q}}z\right)\EXP{-A_n^p \rH{u}^{\frac{pq}{p+q}}(Lpq^{-1}z^{-q}+Lz^p)}\, \d z\\
&\sim&
\sqrt{\frac{2\pi Lp}{p+q}}A_n^{\frac{p}{2}}\rH{u}^{\frac{pq}{2(p+q)}}g\left(A_n \rH{u}^{\frac{q}{p+q}}\right)
\EXP{-D_n \rH{u}^{\frac{pq}{p+q}}},
\end{eqnarray*}
where $D_n=\left(L+Lpq^{-1}\right)A_n^p$ and $A_n=\left(qL_n\right)^{\frac{1}{p+q}}\left(Lp\right)^{-\frac{1}{p+q}}$, and thus the proof is complete.
 \QED

\def\uN{\tilde{u}_n}
\BEL\label{add le3.1}
\cwz{Assume that} the distribution function $F$ of $S$ satisfies
\eqref{eq:SB}, and further \eqref{eq n tail} holds, \cE{then} we have
$$n\sum_{k=1}^{n-1}|\rho(k)|\int_0^{\infty}\int_0^{\infty}
\cH{\exp}\left(-\frac{\left(\uN/s\right)^2+\left(\uN/t\right)^2}{2(1+|\rho(k)|)}\right)\, \d F(s)\, \d F(t)\to
0
$$ as $n \to \infty$, where $\uN=u_n(x)$.
\EEL

\prooflem{add le3.1}
Using similar arguments as in Lemma 4.3.2 in \cite{leadbetter1983extremes},
let $\sigma=\max_{k \ge 1}|\rho(k)|<1$ and $\iota_n=[n^\beta]$, where $\beta$ is any positive constant such that
$\beta < 2(1+\sigma)^{-\frac{p}{2+p}}-1$.
According to \eqref{eq n tail} \wzc{and \eqref{eq:product}} we have
\BQNY
\EXP{-T \uN^{\frac{2p}{2+p}}}\sim \mathbb{Q}\frac{\ccc{g^{-1}\left(Q^{-1}\uN^{\frac{2}{2+p}}\right)}}{n} \qquad \qquad
\uN\sim \left(\frac{\ln n}{T}\right)^{\frac{2+p}{2p}},
\EQNY
\ccc{where $T$ and $Q$ are defined in \eqref{eq T Q}, and $\mathbb{Q}$ is a positive constant which may change from line to line.}

By \eqref{eq:bound SB} \ccc{with $q=2$ and $\cwz{L}_k=1/2(1+|\rho(k)|)$} and split the sum into two parts, i.e.,
\BQNY
&&n\sum_{k=1}^{n-1}|\rho(k)|\int_0^{\infty}\int_0^{\infty}
\cH{\exp}\left(-\frac{\left(\uN/s\right)^2+\left(\uN/t\right)^2}{2(1+|\rho(k)|)}\right)\, \d F(s)\,\d F(t)\\
&\leq&\mathbb{Q}n\sum_{k=1}^{n-1}|\rho(k)|
\uN^{\frac{2p}{2+p}}\ccc{g^{2}\left(\cwz{A_k}\uN^{\frac{2}{2+p}}\right)}
\EXP{-2(1+|\rho(k)|)^{-\frac{p}{2+p}}T\uN^{\frac{2p}{2+p}}}\\
&=&\mathbb{Q}n\left(\sum_{k=1}^{[n^{\beta}]}+\sum_{k=[n^{\beta}]+1}^{n-1}\right)
|\rho(k)|\uN^{\frac{2p}{2+p}}\ccc{g^{2}\left(\cwz{A_k}\uN^{\frac{2}{2+p}}\right)}
\EXP{-2(1+|\rho(k)|)^{-\frac{p}{2+p}}T\uN^{\frac{2p}{2+p}}}.
\EQNY
\ccc{
Since $g(\cdot)$ is ultimately monotone, assume without loss of generality that it is ultimately increasing. By the assumption that  $g(\cdot)$ is a regularly varying function at infinity with index $\alpha$, using Potter bound see e.g., \cite{Res1987}, \cite{deh2006a}, \cite{Faletal2010}  for arbitrary $\varepsilon>0$, \cwz{$k\ge 1$} we have
$$g\left(\cwz{A_k}\uN^{\frac{2}{2+p}}\right)\le g\left(Q^{-1}\uN^{\frac{2}{2+p}}\right)\le \mathbb{Q}\uN^{\frac{2(\alpha+\varepsilon)}{2+p}}$$
}
for all $n$ large. Hence the first part is dominated by
\BQNY
&&\mathbb{Q}nn^{\beta}\uN^{\frac{2p}{2+p}}\ccc{g^2\left(Q^{-1}\uN^{\frac{2}{2+p}}\right)}
\EXP{-2(1+\sigma)^{-\frac{p}{2+p}}T\uN^{\frac{2p}{2+p}}}\\
&=&\mathbb{Q}n^{1+\beta}\uN^{\frac{2p}{2+p}}\ccc{g^2\left(Q^{-1}\uN^{\frac{2}{2+p}}\right)}
\left(\EXP{-T\uN^{\frac{2p}{2+p}}}\right)^{2(1+\sigma)^{-\frac{p}{2+p}}}\\
&\leq&\mathbb{Q}n^{1+\beta}\uN^{\frac{2p}{2+p}}\ccc{g^2\left(Q^{-1}\uN^{\frac{2}{2+p}}\right)}
\left(\frac{\ccc{g^{-1}\left(Q^{-1}\uN^{\frac{2}{2+p}}\right)}}{n}\right)^{2(1+\sigma)^{-\frac{p}{2+p}}}\\
&\leq&\mathbb{Q}n^{1+\beta-2(1+\sigma)^{-\frac{p}{2+p}}} (\ln
n)^{1+\frac{2(\alpha+\varepsilon)}{p}(1-(1+\sigma)^{-\frac{p}{2+p}})}\cH{\to 0}
\EQNY \cH{as $n\to\infty$} since
$1+\beta-2(1+\sigma)^{-\frac{p}{2+p}}<0$.
Next set  $\sigma(l)=\max_{k \ge l}|\rho(k)|<1$. We may further write
\BQNY &&\mathbb{Q}n\sum_{k=\iota_n+1}^{n-1}
|\rho(k)|\uN^{\frac{2p}{2+p}}g^2\left(\cwz{A_k}\uN^{\frac{2}{2+p}}\right)
\EXP{-2(1+|\rho(k)|)^{-\frac{p}{2+p}}T\uN^{\frac{2p}{2+p}}}\\
&\leq&\mathbb{Q}n^2\sigma(\iota_n)
\uN^{\frac{2p}{2+p}}g^2\left(Q^{-1}\uN^{\frac{2}{2+p}}\right)
\EXP{-2(1+\sigma(\iota_n))^{-\frac{p}{2+p}}T\uN^{\frac{2p}{2+p}}}\\
&\leq&\mathbb{Q}n^2\sigma(\iota_n)
\uN^{\frac{2p}{2+p}}g^2\left(Q^{-1}\uN^{\frac{2}{2+p}}\right)
\EXP{-2T\uN^{\frac{2p}{2+p}}}
\EXP{2T\sigma(\iota_n)\uN^{\frac{2p}{2+p}}}\\
&\leq&\mathbb{Q}\sigma(\iota_n)
\uN^{\frac{2p}{2+p}}
\EXP{2T\sigma(\iota_n)\uN^{\frac{2p}{2+p}}}.
\EQNY
Using now \eqref{rho} as $n\to \IF$
\BQNY
\sigma(\iota_n)\uN^{\frac{2p}{2+p}}\sim T^{-1}\sigma(\iota_n)\ln n
\leq T^{-1} \max_{k \ge \iota_n }|\rho(k)| \ln n \to 0
\EQNY
the exponential term above tends to one and the remaining product
tends to zero and thus the proof is complete. \QED

\BEL\label{le ioi}
Let $\{X_n, n \ge 1\}$ be a \eH{ssGs} satisfying \eqref{rho}, and let
$S_n, n \ge 1$ be independent random variables satisfying \eqref{eq:SB} being further independent of $X_n$.
\ccc{Additionally, assume that the survival function of  $Y_n=S_nX_n$ satisfy \eqref{eq n tail}.}
Further if $ 0< \theta <1$ and \ccc{$I_n$ is an interval containing} $k_n \sim \theta n$ members, we have
$$\lim_{n \to \IF}\sup_{x,y \inr} \Abs{\pk{-u_n(y)<m(I_n)\le M(I_n) \le u_n(x)}-\EXP{-\theta (\EXP{-x}+\EXP{-y})}}=0,$$
where $M(I_n)=\EH{\max_{ i \in I_n}}Y_i $ and $m(I_n)=\EH{\min_{  i \in I_n}}Y_i$.
\EEL

\prooflem{le ioi}
Let $Z_n,n\ge 1$ be independent random variables with the same distribution as $X_1$ and define $\mathfrak{M}_n=\max_{1 \le i \le n}S_iZ_i$
and $\mathfrak{m}_n=\min_{1 \le i \le n}S_iZ_i$. For $x,y \inr$, using assumption \eqref{eq n tail}, i.e.,
\BQN\label{eq n*product}
\lim_{n\to \IF}n\pk{S_1Z_1 >u_n(x)}=\EXP{-x}, \quad \lim_{n\to \IF}n\pk{S_1Z_1 \le -u_n(y)}=\EXP{-y}
\EQN
and by Theorem 1.8.2 in \cite{leadbetter1983extremes} we have
\BQN\label{eq iid Mn mn}
\lim_{n \to \IF}\sup_{x,y \inr}\Abs{\pk{-u_n(y)<\mathfrak{m}_n\le\mathfrak{M}_n\le u_n(x)}-\Lambda(x)\Lambda(y)}=0.
\EQN
Further if \eqref{eq:SB} holds, since $S_{n}, n\ge 1$ are independent with common distribution function $F$
by a direct application of \cE{Berman inequality (see \cite{Pit96})}  and Lemma \ref{add le3.1} we obtain
\BQNY
&&\Abs{\pk{-u_n(y)<m_n\le M_n\leq u_n(x)}-
\pk{-u_n(y)<\mathfrak{m}_n\le \mathfrak{M}_n\leq u_n(x)}}\\
&\leq&\int_{[0,\infty]^n}\left|\pk{ \bigcap_{k=1}^n
\left\{-\frac{u_n(y)}{s_k}<X_{k}\leq\frac{u_n(x)}{s_k}\right\}  }
-\pk{  \bigcap_{k=1}^n
\left\{-\frac{u_n(y)}{s_k}<Z_{k}\leq\frac{u_n(x)}{s_k}\right\}}\right|\, \d F(s_1)\cdots \, \d F(s_n)\\
&\leq&
\mathbb{Q}n\sum_{k=1}^{n-1}\int_0^{\infty}\int_0^{\infty}
|\rho(k)|\exp\left(-\frac{(w_n/s)^2+(w_n/t)^2}{2(1+|\rho(k)|)}\right)\, \d F(s) \, \d F(t)\\
&\to & 0, \quad n\to \IF,
\EQNY
where $w_n=\min(|u_n(x)|,|u_n(y)|)$.
Thus by \eqref{eq iid Mn mn} we have
\BQNY\label{eq Mn mn}
\lim_{n \to \IF}\sup_{x,y \inr}\Abs{\pk{-u_n(y)<m_n\le M_n\le u_n(x)}-\Lambda(x)\Lambda(y)}=0.
\EQNY
Now let $v_n=u_{[n/\theta]}$, using \eqref{eq n*product} we get
$$\lim_{n\to \IF}n\pk{S_1Z_1 >v_n(x)}=\theta\EXP{-x}, \quad \lim_{n\to \IF}n\pk{S_1Z_1 \le -v_n(y)}=\theta\EXP{-y},$$
hence
\BQN\label{M vn}
\lim_{n \to \IF}\sup_{x,y \inr}\Abs{\pk{-v_n(y)<m_n\le M_{n} \le v_n(x)}-\EXP{-\theta (\EXP{-x}+\EXP{-y})}}=0.
\EQN
Since $S_{n},n\ge 1$ are independent and have a common distribution function $F$, by the stationarity of $\{X_n, n \ge 1\}$
\BQNY
\pk{-u_n(y)<m(I_n)\le M(I_n) \le u_n(x)}&=&\pk{\bigcap_{i \in I_n}\{-u_n(y)<S_iX_i \le u_n(x)\}}\\
&=&\int_{(0, \IF)^{k_n}}\pk{\bigcap_{i=1}^{k_n}\left\{-\frac{u_n(y)}{s_i} <X_i \le \frac{u_n(x)}{s_i}\right\}}\, \d F(s_1)\cdots \d F(s_{k_n})\\
&=&\pk{-u_n(y)<m_{k_n}\le M_{k_n} \le u_n(x)}.
\EQNY
Hence, replacing $n$ by $k_n$ in \eqref{M vn} establishes the claim .
\QED

\BRM\label{rm ioi Mn}
Under the conditions of Lemma \ref{le ioi}, we have
$$\lim_{n \to \IF}\sup_{x \inr}\Abs{\pk{ M(I_n) \le u_n(x)}-\EXP{-\theta \EXP{-x}}}=0.$$
\ERM

\BEL\label{le Mn mn In}
Under the assumptions of Lemma \ref{le ioi}, and
let $I_1,I_2,\ldots,I_l$ (with $l$ a fixed number) be disjoint subintervals of $\{1,2,\ldots,n\}$ such that $I_i$ has
$k_{n,i}\sim \theta_i n$
elements, where $\theta_i$ are fixed positive constants with $\theta:=\sum^l_{i=1}\theta_i \le 1$. Then we have
$$\pk{\bigcap_{i=1}^l\{-u_n(y)< m(I_i) \le M(I_i)\le u_n(x)\}}-\prod_{i=1}^l\pk{-u_n(y)< m(I_i) \le M(I_i)\le u_n(x)}\to 0$$
as $n \to \IF$.
\EEL

\prooflem{le Mn mn In}
Since $\{X_n, n \ge 1\}$ is a stationary random sequence, using Berman's inequality 
 and Lemma \ref{add le3.1}, we have
\BQNY
&&\left|\pk{\bigcap_{i=1}^l\{-u_n(y)< m(I_i) \le M(I_i)\le u_n(x)\}}-\prod_{i=1}^l\pk{-u_n(y)< m(I_i) \le M(I_i)\le u_n(x)}\right|\\
&=&\left|\pk{\bigcap_{i=1}^l\bigcap_{j \in I_i}\{-u_n(y)< S_j X_j\le u_n(x)\}}-\prod_{i=1}^l\pk{\bigcap_{j \in I_i}\{-u_n(y)< S_j X_j\le u_n(x)\}}\right|\\
&\le&\int_{(0,\IF)^{\hat{\theta}_l}}\left|\pk{\bigcap_{i=1}^l \hat{A}_i}-\prod_{i=1}^l\pk{\hat{A}_i}\right|\,\d F(s_1)\cdots \d F(s_{\hat{\theta}_l})\\
&\le& \hat{\theta}_l \sum_{k=1}^{\hat{\theta}_l}\int_0^{\IF}\int_0^{\IF}|\rho(k)|\exp\left(-\frac{(w_n/s)^2+(w_n/t)^2}{2(1+|\rho(k)|)}\right) \,\d F(s) \d F(t)\\
&\to& 0, \quad n \to \IF,
\EQNY
where $\hat{A}_i=\bigcap_{j=\hat{\theta}_{i-1}+1}^{\hat{\theta}_i}\left\{-\frac{u_n(y)}{s_j}<  X_j\le \frac{u_n(x)}{s_j}\right\}$ with $$\hat{\theta}_i=\sum_{j=1}^i [\theta_j n], \quad \hat{\theta}_{0}=0, \quad w_n=\min(|u_n(x)|,|u_n(y)|),$$
hence the proof is complete. \QED

\BRM\label{rm Mn In}
Under the condition of Lemma \ref{le Mn mn In}, we have
$$\limit{n} \ABs{ \pk{\bigcap_{i=1}^l\{ M(I_i)\le u_n(x)\}}-\prod_{i=1}^l\pk{ M(I_i)\le u_n(x)}}= 0.$$
\ERM

\prooftheo{th N}
According to  \cite{Win96}, first we need to prove
that the marginal point processes of $N_{n,d}$ converge weakly to a Poisson process $N_{d}$ with intensity $\EXP{-x_d},d=1,2$.
By Theorem A.1 in \cite{leadbetter1983extremes} for $N_{n,1}(B_1,x_1)$, it is sufficient to show that \cE{as $n\to \IF$}
\BQNY
&(P_1).&
\EE{ N_{n,1}((s,t],x_1)} \to \EE{ N_1((s,t],x_1)}=(t-s)\EXP{-x_1}, \, 0<s<t\le 1;\\
&(P_2).& \pk{\bigcap_{i=1}^k\{N_{n,1}((s_i, t_i],x_1)=0\}} \to
\pk{\bigcap_{i=1}^k\{N_1((s_i, t_i],x_1)=0\}} =
\EXP{-\sum_{i=1}^k(t_i-s_i)\EXP{-x_1}},
\EQNY
where $0<s_1<t_1 \le s_2 < t_2  \le \cdots \le s_k <t_k\le 1$.\\
We have
\BQNY \EE{ N_{n,1}((s,t],x_1)}&=& \EE{ \sum_{i/n
\in (s,t]} \I(S_iX_i >u_n(x_1))}
=\sum_{i/n \in (s,t]} \pk{S_iX_i >u_n(x_1)}\\
&\to& (t-s)\EXP{-x_1}=\EE{N_1((s,t],x_1)}
\EQNY
as $n \to \IF$, where the above convergence follows from \eqref{eq n tail}.\\
In order to show $(P_2)$ note first that for $0<s<t \le 1$
$$\pk{N_{n,1}((s,t],x_1)=0}=\pk{ M(I_n)\le u_n(x_1)},$$
where $I_n=\{[sn]+1,\ldots,[tn]\}$.
Further, $I_n$ contains $k_n$ integers with $k_n=[tn]-[sn] \sim (t-s)n$ as $n \to \IF$. Thus,
 in view of Remark \ref{rm ioi Mn} with $\theta=t-s<1$ we have
 \BQN\label{eq Nn
In2} \pk{N_{n,1}((s,t],x_1)=0} \to \EXP{-(t-s)\EXP{-x_1}} \quad
\mbox{as } n \to \IF.
\EQN
Next, let $E_i$ be the set of integers
$\{[s_{i}n]+1, \ldots, [t_{i}n]\}$ with $0<s_1<t_1 \le s_2 < t_2 \le
\cdots \le s_k <t_k\le 1$, then we have
\BQNY
\lefteqn{\pk{\bigcap_{i=1}^k\{N_{n,1}((s_i,t_i],x_1)=0\}}}\\
&=&\pk{\bigcap_{i=1}^k\{ M(E_i)\le u_n(x_1)\}}\\
&=&\prod_{i=1}^k\pk{N_{n,1}((s_i,t_i],x_1)=0}+\left(\pk{\bigcap_{i=1}^k\{
M(E_i)\le u_n(x_1)\}}-\prod_{i=1}^k\pk{ M(E_i)\le u_n(x_1)}\right).
\EQNY
Using \eqref{eq Nn In2}, the first term converges to
$\EXP{-\sum_{i=1}^k(t_i-s_i)\EXP{-x_1}}$ as $n \to \IF$.
 By Remark \ref{rm Mn In} the modulus of the remaining difference of terms tend to 0.
Consequently,  $N_{n,1}$ converge weakly to a Poisson process $N_1$ with intensity $\EXP{-x_1}$.
Since $Y_i \EQD -Y_i$, $N_{n,2}$ also converge weakly to a Poisson process $N_2$ with intensity $\EXP{-x_2}$.

\peng{Now define the avoidance function of $\mathbf{N}_n$ as \BQNY
F_{\mathbf{{N}}_n}(\mathbf{B})=\pk{N_{n,1}(B_1,x_1)=0,
N_{n,2}(B_2,x_2)=0}, \EQNY where $B_{1}$ and $B_{2}$ are defined
below. To get the main result, it suffices to prove that \BQNY
\lim_{n \to \IF}F_{\mathbf{N}_n}(\mathbf{B}) \mbox{ exists for all }
\mathbf{B}=\bigcup_{d=1}^2\bigcup_{j=1}^m(B_{dj}\times \{d\}), \EQNY
for arbitrary positive integers $m$, where $B_{dj}=(s_{dj},t_{dj}]$,
$0< s_{d1} < t_{d1} \le s_{d2} <t_{d2} \le \ldots \le s_{dm} <t_{dm}
\le 1$, and $B_1=\bigcup_{j=1}^m B_{1j}$, $B_2=\bigcup_{j=1}^m
B_{2j}$. We will show that \[\lim_{n\to
\IF}F_{\mathbf{N}_n}(\mathbf{B})=\EXP{-m(B_{1})\EXP{-x_{1}}
-m(B_{2})\EXP{-x_{2}}}.\] For simplicity we only consider the case
$B_{1}\subset B_{2}$;  other cases are similar. First consider the
case $n(B_{2}\setminus B_{1})=o(n)$, i.e., $m(B_1)=m(B_2)$.
Obviously,
\begin{eqnarray*}0&\le&\pk{-u_{n}(x_2)<Y_k \le u_n(x_1), k/n \in B_1}-\pk{Y_k \le u_n(x_1), k/n \in B_1; -Y_l \le
u_n(x_2), l/n \in
B_2}\\
&\le&\sum_{l:l/n\in B_{2}\setminus B_{1}}\pk{-Y_{l}>u_{n}(x_2)}\to0
\end{eqnarray*}
as $n\to\infty$.  Consequently,
by Lemma \ref{le ioi} and \ref{le Mn mn In}, we have
\begin{eqnarray*}
&&\lim_{n\to\infty}\pk{Y_k \le u_n(x_1), k/n \in B_1; -Y_l \le
u_n(x_2), l/n \in B_2}\\
&=& \lim_{n\to\infty}\pk{-u_{n}(x_2)<Y_k \le u_n(x_1), k/n \in
B_1}\\
&=&\prod_{j=1}^m \EXP{-(t_{1j}-s_{1j})\EXP{-x_1}}\prod_{j=1}^m
\EXP{-(t_{1j}-s_{1j})\EXP{-x_2}}\\
&=&\EXP{-m(B_{1})\EXP{-x_1}-m(B_{2})\EXP{-x_2}}.
\end{eqnarray*}
It suffices to prove the case of $n(B_{2}\setminus B_{1})=O(n)$.
Noting that any $z>0$, \BQNY &&\mathbb{P}\left(-u_n(x_2) < Y_k \le
u_n(x_1), k/n \in B_1;-u_n(x_2) < Y_i \le u_n(z), i/n \in B_2
\setminus B_1\right)\\&\le&\pk{Y_k \le u_n(x_1), k/n \in B_1; -Y_l
\le u_n(x_2), l/n \in B_2}\\ &\le&
\mathbb{P}\left(-u_n(x_2) < Y_k \le u_n(x_1), k/n \in B_1;-u_n(x_2) < Y_i \le u_n(z), i/n \in B_2 \setminus B_1\right)\\
&&+\mathbb{P}\left(\max(Y_{i},i/n \in B_{2}\setminus B_{1})>u_{n}(z)\right)\\
&=&\mathbb{P}\left(-u_n(x_2) < Y_k \le u_n(x_1), k/n \in
B_1;-u_n(x_2) < Y_i \le u_n(z), i/n \in B_2 \setminus
B_1\right)\\
&&+\left(1-\mathbb{P}\left(\max(Y_{i},i/n \in B_{2}\setminus
B_{1})\le u_{n}(z)\right)\right). \EQNY
\cE{Applying} Lemma \ref{le ioi}
and \ref{le Mn mn In} once again, we obtain
\begin{eqnarray*}
&&\EXP{-m(B_{1})\left(\EXP{-x_{1}}+\EXP{-x_{2}}\right)}\EXP{-m(B_{2}\setminus B_{1})\left(\EXP{-z}+\EXP{-x_{2}}\right)}\\
&\le &\liminf_{n\to\infty}\pk{Y_k \le u_n(x_1), k/n \in B_1; -Y_l
\le u_n(x_2), l/n \in B_2}\\
&\le& \limsup_{n\to\infty}\pk{Y_k \le u_n(x_1), k/n \in B_1; -Y_l
\le
u_n(x_2), l/n \in B_2}\\
&\le&
\EXP{-m(B_{1})\left(\EXP{-x_{1}}+\EXP{-x_{2}}\right)}\EXP{-m(B_{2}\setminus B_{1})\left(\EXP{-z}+\EXP{-x_{2}}\right)}\\
&&\qquad+\left(1-\EXP{-m(B_{2}\setminus B_{1})\EXP{-z}}\right).
\end{eqnarray*}
\cE{Hence}, letting $z\to\infty$ we have \[\lim_{n\to\infty}\pk{Y_k \le
u_n(x_1), k/n \in B_1; -Y_l \le u_n(x_2), l/n \in
B_2}=\EXP{-m(B_{1})\EXP{-x_{1}}-m(B_{2})\EXP{-x_{2}}}\] }
\cE{establishes the proof}. \QED

\proofkorr{co M m}
Notice that
\BQNY
\pk{M_n^{(k)} \le u_n(x), m_n^{(l)} > -u_n(y)}=\pk{N_{n,1}((0,1],x) \le k-1, N_{n,2}((0,1],y) \le l-1}.
\EQNY
Hence the proof follows by an immediate application of Theorem \ref{th N}. \QED

\prooftheo{th N bound} According to Lemma 3.3 of \cite{HashorvaWengB}, i.e., under the condition \eqref{eq rho},
$$\lim_{n \to \IF}n\sum_{k=1}^{n-1}|\rho(k)|\int_0^{1}\int_0^{1}
\cH{\exp}\left(-\frac{\left(u_n(x)/s\right)^2+\left(u_n(x)/t\right)^2}{2(1+|\rho(k)|)}\right)\, \d F(s)\, \d F(t)=
0$$
is valid, then
we have Lemma \ref{le ioi} and \ref{le Mn mn In} also hold for $S_n$ satisfy \eqref{eq:bound S}. Hence
using the similar argument as Theorem \ref{th N}, the desired result obtained.
\QED

\vspace{0.5in}

\textbf{Acknowledgments.} Z. Weng has been partially supported by the Swiss
National Science Foundation grant 200021-134785 and by
the project RARE -318984 (a Marie Curie IRSES Fellowship within the 7th European Community
Framework Programme);
 Z. Peng  has been supported by the National Natural Science Foundation of China under grant
11171275 and the Natural Science Foundation Project of CQ under
cstc2012jjA00029.

\bibliographystyle{plain}
\bibliography{Gauss}
\end{document}